\documentclass[titlepage,twoside,12pt]{article}
\usepackage{amssymb}
\usepackage{amsfonts}
\begin{document}

{\Large \bf What is wrong with the Lax-Richtmyer \\ \\ Fundamental Theorem of Linear \\ \\
Numerical Analysis ?} \\

Elemer E Rosinger \\
Department of Mathematics \\
and Applied Mathematics \\
University of Pretoria, Pretoria \\
0002 South Africa \\
eerosinger@hotmail.com \\ \\

{\bf Abstract} \\

We show that the celebrated 1956 Lax-Richtmyer linear theorem in Numerical Analysis - often
called the {\it Fundamental Theorem of Numerical Analysis} - is in fact wrong. Here "wrong"
does not mean that its statement is false mathematically, but that it has a limited practical
relevance as it misrepresents what actually goes on in the numerical analysis of partial
differential equations. Namely, the assumptions used in that theorem are excessive to the
extent of being unrealistic from practical point of view. The two facts which the mentioned
theorem gets wrong from practical point of view are : \\
- the relationship between the convergence and stability of numerical methods for linear
PDEs, \\
- the effect of the propagation of round-off errors in such numerical methods. \\
The mentioned theorem leads to a result for PDEs which is unrealistically better than the well
known best possible similar result in the numerical analysis of ODEs. Strangely enough, this
fact seems not to be known well enough in the literature. Once one becomes aware of the above,
new avenues of both practical and theoretical interest can open up in the numerical analysis
of PDEs. \\ \\

\newpage

{\bf 1. Towards a correct relationship between stability and \\
\hspace*{0.5cm} convergence} \\

It has been shown that in practically relevant situations the converse implication "convergent
$\Longrightarrow$ stable" in the Lax-Richtmyer theorem may {\it fail} to hold, see Rosinger
[1-8], Rosinger \& van Niekerk [1,2], Oberguggenberger \& Wang. Thus there need {\it not}
always be an equivalence between the {\it convergence} and {\it stability} of a numerical
scheme. It may therefore happen that convergence is a {\it weaker} property than stability,
which means that we may have convergent numerical schemes which nevertheless {\it fail} to be
stable as well. \\

In this way, what has become a kind of "UNIVERSALLY RECITED MANTRA" in the numerical analysis
of linear partial differential equations, namely that \\

( ? )~~~ {\it "stability and convergence are equivalent"}  \\

for linear numerical methods approximating such equations, does in fact {\it lack} a valid
enough practical reality, and can be replaced with the far more convenient fact, according to
which \\

( ! )~~~ {\it "convergence need not always imply stability"}. \\

Needless to say, the practical interest in such a possibility is significant, as it can {\it
enlarge} the class of convergent numerical schemes {\it beyond} those which are stable.
Examples in this regard are mentioned in Rosinger [1-8], Rosinger \& van Niekerk [1,2],
Oberguggenberger \& Wang. \\

By the way, the well known {\it necessary} condition for stability, given by von Neumann prior
to the Lax-Richtmyer theorem, and which does {\it not} require a Banach space setup, can be
seen as a further indication of the rather involved relationship between stability and
convergence. \\

A yet more important point about the above mantra is the following. Even if it were true in
the linear case - which in fact is the only case addressed by the Lax-Richtmyer theorem - it
would still lack relevance in most of the cases when exact solutions of {\it nonlinear} PDEs
are approximated by respective {\it nonlinear} numerical methods. Indeed, as is well known,
Kreiss, Stetter, a local linearized stability analysis of nonlinear PDEs and of their
nonlinear numerical methods need {\it not} in general lead either to necessary, or sufficient
convergence conditions. \\

{\bf 2. Questions about the implication "convergent} $\Longrightarrow$ {\bf stable"} \\

Briefly, the Lax-Richtmyer theorem, see below, states the {\it equivalence} between the {\it
convergence} and {\it stability} of a linear numerical scheme which is consistent with a well
posed linear PDE, see Lax \& Richtmyer, Richtmyer, Richtmyer \& Morton, as well as a fully
detailed analysis and presentation in Rosinger [2, pp. 1-14]. \\

The important fact to note is the following. \\

The proof of the implication "stable $\Longrightarrow$ convergent" is trivial, and certainly,
it does in {\it no way} require the completeness of space in which it happens. \\
Therefore, the crux of the Lax-Richtmyer theorem is solely in the proof of the converse
implication, namely, "convergent $\Longrightarrow$ stable". \\

That converse implication "convergent $\Longrightarrow$ stable", however, is proved based on
the celebrated Principle of Uniform Boundedness of Linear Operators in Banach spaces. And as
is well known, see Appendix, that property of uniform boundedness does {\it not} necessarily
hold in normed spaces which are not complete, thus fail to be Banach spaces. \\

It is precisely here, with the assumptions which are made in order to secure a Banach space
framework, that the Lax-Richtmyer theorem goes {\it twice} wrong from practical point of view.
Namely, it goes wrong both with respect to the relationship between stability and convergence,
as well as regarding the treatment of the essentially nonlinear phenomenon of the propagation
of round-off errors. \\

Furthermore, the proof of the implication "convergent $\Longrightarrow$ stable" is essentially {\it linear}, as it
makes use of the mentioned linear principle, as well as of a linear concept of stability. This makes the extension of
that implication to the fully nonlinear case extremely difficult. \\

{\bf 3. Is completeness an appropriate requirement ?} \\

The numerical analysis of a given PDE does typically assume the a priori knowledge of the
existence of certain exact solutions of that equation. After all, in case an exact solution
does not exist, it is of course nonsensical to try to approximate it numerically. Thus, if and
when the existence of an exact solution is known a priori, then the aim of numerical analysis
is to construct numerical solutions approximating one or another of such exact solutions. In
this way, we are given an exact solution $U$ and construct a sequence, say, $U_{\Delta t}$,
with $\Delta t > 0$, of numerical solutions. Thus the problem is whether or not we have the
hoped for convergence property \\

( * ) ~~~ $ \lim_{\Delta t \to 0}~ U_{\Delta t} ~=~ U $ \\

where the limit holds in some appropriate sense. \\

Suppose now, as usual, that both the exact solution $U$ and the numerical solutions
$U_{\Delta t}$ belong to a certain normed space $( X, ||~||~ )$. \\

A {\it crucial observation} here is the following one. And it is {\it missed} by the
Lax-Richtmyer theorem. \\

Clearly, in order to establish whether the above convergence property ( * ) does, or for that
matter, does not hold, one does {\it not at all} need to assume that the respective normed
space $( X, ||~||~ )$ is {\it complete}. Indeed, we have started by assuming that the exact
solution $U$ exists, thus the hoped for limit value in (*) exists. Furthermore, the terms
$U_{\Delta t}$ of the sequence in (*) also exist, being the constructed numerical solutions.
Finally, the normed space $X$ is supposed to be chosen in such a way that both the exact
solution $U$ and its numerical approximations $U_{\Delta t}$ do belong to it. And then the
only problem is whether the constructed numerical solution $U_{\Delta t}$ does indeed happen
to converge to the exact solution $U$. \\

Furthermore, often, when for instance the exact solution $U$ is a classical solution of the
PDE considered, one can choose the normed space $( X, ||~||~ )$ as constituted from
sufficiently {\it smooth} functions, since the numerical solutions $U_{\Delta t}$ are
typically defined at discrete points, thus they can be extrapolated to functions of required
smoothness. \\

In this way, there does {\it not} appear to be any practical reason whatsoever why the normed
space $( X, ||~||~ )$ in which the convergence property (*) is to be established should be
complete. \\

The alleged reason why nevertheless the completeness of the normed space $( X, ||~||~ )$ is
requested appears to be the claim that it is needed in order to handle the effect of round-off
errors as well. Indeed, as it stands, the Lax-Richtmyer theorem is only supposed to deal with
the effects of the propagation of truncation errors, since it does not anywhere mention
directly round-off errors. \\
However, as seen in Rosinger [5], see also Rosinger [2-4,7], this claim that the completeness
of $X$ will give the opportunity to deal as well with the effect of the propagation of
round-off errors is simply {\it unrealistic} from the point of view of the way round-off
errors actually propagate in the computations involved. In particular, this claim leads to the
{\it paradox} that one obtains a result regarding the effect of round-off errors in the
numerical solution of PDEs which is strictly better than the well known best possible
corresponding result in the case of the numerical solution of ODEs. \\

Obviously, in the many usual cases when one approximates classical solutions $U$, one can
choose the normed space $( X, ||~||~ )$ made up of sufficiently smooth functions. But then,
the completeness requirement in the Lax-Richtmyer theorem obliges one to consider its
completion $( X^{\#}, ||~||~ )$. And typically, $X^{\#}$ will be a much larger space,
containing a considerable amount of non-smooth functions. \\

Two important points should be noted here. \\

First, within this larger and completed space $X^{\#}$, the original {\it convergence} problem
( * ) will remain precisely the {\it same}. Indeed, the constructed sequence of numerical
solutions $U_{\Delta t} \in X$ converges to the existing exact solution $U \in X$ in the space
$X$, if and only if it converges to $U$ in the space $X^{\#}$. \\

On the other hand, the {\it stability} property of the respective numerical method may turn
out to lead to a more {\it stringent} condition in the larger space $X^{\#}$, than in the
original smaller space $X$. \\

This is indeed of one the issues related to the assumption of completeness, an assumption
which is essential in the particular method of proof of the implication "convergent
$\Longrightarrow$ stable" in the Lax-Richtmyer theorem. \\
Otherwise, one simply notes that, in general, the completeness condition does {\it not}
necessarily belong to the problem of establishing the convergence property (*). \\

{\bf 4. Compactness or Boundedness ?} \\

The above {\it convergence} relation ( * ), whenever it holds, clearly implies that the
subset \\

$~~~ \{~ U_{\Delta t} ~~|~~ \Delta t > 0 ~\} \cup \{~ U ~\} \subset X $ \\

is {\it compact} in $X$, regardless of $X$ being complete or not. And let us recall that all
the elements of this subset are supposed to exist. Indeed, the exact solution $U$ of the PDE
under consideration exists, otherwise the problem of its numerical approximation would be
vacuous. Further, the approximating numerical solutions $U_{\Delta t}$ are effectively
constructed by the numerical method employed. \\

On the other hand, the condition of {\it stability} of the numerical methods used in the
Lax-Richtmyer theorem, see (5.18) below, is given in terms of {\it boundedness}, and as is
well known, boundedness does {\it not} imply compactness in infinite dimensional normed
spaces. \\

This {\it discrepancy} between the association of convergence with compactness, and on the
other hand, of stability with boundedness was first pointed out and dealt with in Rosinger [1],
where with an appropriate compactness based definition of stability, a general {\it nonlinear
equivalence} result was given between convergence and stability. \\

It should be mentioned here that the above arguments related to stability, convergence,
completeness, compactness and boundedness were, back in the early summer of 1979, personally
communicated by the author to P D Lax, at a conference at the Tel Aviv University, in
Israel. \\

{\bf 5. Some details of the Lax-Richtmyer theorem} \\

Let us now, for convenience, recall the Lax-Richtmyer theorem as given in its original
formulation, see Lax \& Richtmyer, Richtmyer, Richtmyer \& Morton, Rosinger [2, pp. 1-14]. We
consider a {\it linear evolution} type PDE \\

(5.1)~~~ $ d/dt~ U ( t ) ~=~ A ( U ( t ) ),~~ t \in [ 0, T ] $ \\

with the initial value \\

(5.2)~~~ $ U ( 0 ) ~=~ u $ \\

where $A : D \subseteq X \longrightarrow X$ is a linear operator defined on the subspace $D$
of the Banach space $X$, $u \in D$, while $U : [ 0, T ] \longrightarrow D$ is the sought after
solution. Since we deal with an evolution PDE, the operator A is in fact a linear partial
differential operator in some space variable $x \in \mathbb{R}^n$. \\
Further, one can assume that, when given, linear homogenous boundary conditions have already
been incorporated in the definition of $D$. \\
Typically, one can also assume that $D$ is dense in $X$ and we have satisfied the following
{\it exact solution} property \\

(5.3)~~~ $ \begin{array}{l}
               \forall~~~ u \in D ~: \\ \\
               \exists~~~ U : [ 0, T ] ~\longrightarrow~ X ~: \\ \\
               ~~* )~~~ \lim_{\Delta t \to 0}~ ||~
                           ( U ( t + \Delta t ) - U ( t ) ) / \Delta t - A ( U ( t ) ) ~||
               ~=~ 0,~~~ t \in [ 0, T ] \\ \\
               ~** )~~~ U ( 0 ) ~=~ u
            \end{array} $ \\

Given now {\it time}, respectively {\it space} increments $\Delta t \in ( 0, \infty )$ and
$\Delta x \in ( 0, \infty )^n$, we construct a {\it finite difference} method \\

(5.4)~~~ $ C_{\Delta t,~ \Delta x} : X ~\longrightarrow~ X $ \\

which we assume to be a continuous linear mapping. \\

The {\it numerical analysis problem} we face in the above terms is to {\it characterize} the
relations \\

(5.5)~~~ $ \Delta x ~=~ \alpha ( \Delta t ) $ \\

where the mapping $\alpha : ( 0, \infty ) \longrightarrow ( 0, \infty )^n$ is such that
$\lim_{\Delta t \to 0}~ \alpha ( \Delta t ) = 0 \in \mathbb{R}^n$, and the {\it convergence}
property holds \\

(5.6)~~~ $ \lim_{\Delta t \to 0,~ n \to \infty,~ n \Delta t \to t}~ ||~ U ( t ) -
                 C^n_{\Delta t,~ \alpha ( \Delta t )}~ u ~||~=~ 0 $ \\

uniformly for $t \in [ 0, T ]$, for every $u \in D$, where $U$ corresponds to $u$ according to
(5.3). \\

As is well known, in general, this is not a trivial problem. In Courant et.al., it was shown
for the first time that one {\it cannot} in general expect instead of (5.6) the stronger
convergence property \\

(5.7)~~~ $ \lim_{\Delta t \to 0,~ n \to \infty,~ n \Delta t \to t,~ \Delta x \to 0}~
                    ||~ U ( t ) - C^n_{\Delta t,~ \Delta x}~ u ~||~=~ 0 $ \\

to hold uniformly for $t \in [ 0, T ]$. \\

Property (5.7), in which $\Delta t$ and $\Delta x$ can simultaneously and independently tend
to 0, is called {\it unconditional stability}. On the other hand, property (5.6), in which the
relation (5.5) ties $\Delta x$ to $\Delta t$ when they both tend to 0, is called {\it
conditional stability}. \\

Obviously, in the case of conditional stability, one is interested in numerical methods (5.4)
in which $\alpha$ tend to 0 as fast as possible, when $\Delta t$ tends to 0. Indeed, in such a
situation one can obtain a good space accuracy without increasing too much the computation
time. \\

As a simple and immediate illustration, let us consider the initial value problem for the heat
equation \\

~~~$ U_t ~=~ U_{xx},~~~ t \in [ 0, \infty ),~ x \in \mathbb{R} $ \\

~~~$ U ( 0, x ) ~=~ u ( x ),~~~ x \in \mathbb{R} $ \\

In this case we can take $( X, ||~|| ~) = {\cal L}^\infty ( \mathbb{R} )$ and $A = \partial^2 /
\partial x^2$, while \\

~~~$ D ~=~ \{~ u \in {\cal L}^\infty ( \mathbb{R} )\cap {\cal C}^2 ( \mathbb{R} ) ~~|~~
                                          A u \in {\cal L}^\infty ( \mathbb{R} ) ~\} $ \\

and as is well known, John, the exact solution property (5.3) holds. \\
A simple numerical method for this heat equation is given by \\

~~~$ ( C_{\Delta t,~ \Delta x}~ u ) ( x ) ~=~ u ( x ) + ( u ( x + \Delta x ) -
              2 u ( x ) + u ( x - \Delta x ) )
                       \Delta t / \Delta x^2 $ \\

for $u \in X,~ x \in \mathbb{R},~ \Delta t,~ \Delta x > 0$. Also as is well known, Richtmyer,
this explicit numerical method will {\it not} have the convergence property given by the
unconditional stability (5.7), while the weaker convergence property (5.6), called conditional
stability, will hold, if and only if \\

~~~$ 2 \Delta t \leq \Delta x^2 $ \\

which in terms of (5.5) can be written as \\

~~~$ \alpha ( \Delta t ) \geq \sqrt ( 2 \Delta t ) $ \\

This obviously means that $\alpha \to 0$ rather slowly, when $\Delta t \to 0$, which is
inconvenient, since a small $\Delta x$ will impose the use of a quadratically smaller
$\Delta t$, leading thus to an increased number of time iterations. \\

Returning to the general case, in view of the exact solution property (5.3), we can define the
family of linear mappings \\

(5.8)~~~ $ E_0 ( t ) : D ~\longrightarrow~ X,~~~ t \in [ 0, T ] $ \\

by \\

(5.8.1)~~~ $ E_0 ( t ) ( u ) ~=~ U ( t ),~~~ t \in [ 0, T ] $ \\

The initial value problem (5.1), (5.2) is called {\it properly posed}, if and only if the
family of linear mappings (5.8) is uniformly bounded, that is, for a certain $K > 0$, we
have \\

(5.9)~~~ $ ||~ E_0 ( t ) ~|| \leq K,~~~ t \in [ 0, T ] $ \\

As is well known, since $D$ is a dense subspace in $X$, one can extend by continuity the
family of linear mappings (5.8) to a unique family of linear mappings \\

(5.10)~~~ $ E ( t ) : X ~\longrightarrow~ X,~~~ t \in [ 0, T ] $ \\

with the same uniform bound, namely \\

(5.11)~~~ $ ||~ E ( t ) ~|| \leq K,~~~ t \in [ 0, T ] $ \\

In addition, we shall have the {\it semigroup} property \\

(5.12)~~~ $ \begin{array}{l}
                      E ( t ) E ( s ) ~=~ E ( t + s ),~~~
                                        t, s \in [ 0, T ],~ t + s \leq T \\ \\
                      E ( 0 ) ~=~ \mbox{id}_X \\ \\
                      \lim_{\Delta t \to 0}~ ||~ E ( u ) - u ~|| ~=~ 0,~~~ u \in X
            \end{array} $ \\ \\

We can note that the semigroup property (5.12) leads to a further extension, this time of the
linear mappings (5.10), namely \\

(5.13)~~~ $ E ( t ) : X ~\longrightarrow~ X,~~~ t \in [ 0, \infty ) $ \\

where \\

(5.13.1)~~~ $ E ( t ) ~=~ E ( t - [ t / T ] T ) E ( T ) ^ {[ t / T ]},~~~
                                                      t \in [ T, \infty ) $ \\

with $[ t / T]$ denoting the largest integer which is smaller than, or equal to $t / T$. In
this case, instead of the corresponding above relations, we shall have \\

(5.14)~~~ $ \begin{array}{l}
                        E ( t ) E ( s ) ~=~ E ( t + s ),~~~ t, s \in [ 0, \infty ) \\ \\
                        ||~ E ( t ) ~|| \leq K^{1 + [ t / T ]},~~~ t \in [ 0, \infty )
             \end{array} $ \\ \\

Returning now to the numerical method (5.4), (5.5), we shall consider as our {\it finite
difference} scheme the family of continuous linear mappings \\

(5.15)~~~ $ C_{\Delta t} ~=~ C_{\Delta t,~ \alpha ( \Delta t )} : X ~\longrightarrow~ X $ \\

Here it is important to note that, typically, this family of continuous linear mappings
$C_{\Delta t}$, with $\Delta t > 0$, is {\it not} uniformly bounded for small $\Delta t$. \\

Now in view of (5.5), (5.6), the finite difference scheme (5.15) is called {\it convergent} to
the semigroup (5.10) - (5.12) on the time interval $[ 0, T ]$, where $T > 0$ is given, if and
only if \\

(5.16)~~~ $ \begin{array}{l}
                        \forall~~~ u \in X,~ \epsilon > 0 ~: \\ \\
                        \exists~~~ \delta > 0 ~: \\ \\
                        \forall~~~ t \in [ 0, T ],~ \Delta t > 0,~ n \in
                                      {\bf N},~ n \Delta t \leq T ~: \\ \\
                        ~~~~ \Delta t,~ |~ t - n \Delta t ~| \leq \delta
                                           ~~~\Longrightarrow ~~~
                                    ||~ E ( t ) - C^n_{\Delta t} u ~|| \leq \epsilon
             \end{array} $ \\ \\

Further, the finite difference scheme (5.15) is called {\it consistent} with the initial value
problem (5.1), (5.2) on the same time interval $[ 0, T ]$, if and only if \\

(5.17)~~~ $ \begin{array}{l}
                        \forall~~~ u \in X,~ \epsilon > 0 ~: \\ \\
                        \exists~~~ \theta > 0 ~: \\ \\
                        \forall~~~ t \in [ 0, T ],~ \Delta t > 0 ~: \\ \\
                        ~~~~ \Delta t \leq \theta ~~~\Longrightarrow ~~~
                        ||~ C_{\Delta t}~ E ( t ) - E ( t + \Delta t ) u ~|| \leq \epsilon
             \end{array} $ \\ \\

Finally, the finite difference scheme (5.15) is called {\it stable} on the time interval
$[ 0, T ]$, if an only if \\

(5.18)~~~ $ \begin{array}{l}
                       \exists~~~ L > 0 ~: \\ \\
                       \forall~~~ \Delta t > 0,~ n \in {\bf N},~ n \Delta t \leq T ~: \\ \\
                       ~~~~ ||~ C^n_{\Delta t} ~|| \leq L
             \end{array} $ \\ \\

With the above, we have the so called {\it Fundamental Theorem of Linear Numerical Analysis} \\

{\bf Theorem ( Lax-Richtmyer, 1956)} \\

Given a properly posed semigroup (5.10) - (5.12) and a finite difference scheme (5.15) which
is consistent with it, then the finite difference scheme is convergent to the semigroup, if
and only if it is stable. \\

{\bf Remark} \\

The practical interest in the above type of result is in the following. The {\it consistency}
of a finite difference scheme with a semigroup generated by an initial value problem (5.1),
(5.2) is typically easy to establish with the use of a finite Taylor series argument, in case
we deal with smooth enough, or classical solutions. Also, what is practically particularly
important, the consistency property can be established {\it without} the effective knowledge
of any specific exact solution of the initial value problem, and only based on the knowledge
of the regularity of such solutions, that is, the existence of smooth enough, or classical
exact solutions. \\
The {\it convergence} property of such a finite difference scheme is, of course, the main and
nontrivial issue, and just like the consistency property, it is a {\it relational} property,
since it involves the semigroup, or the initial value problems as well. Furthermore, here the
fact that, typically, the exact solution is only known to exist, but it is not known
effectively - this being the very reason for using numerical analysis - makes it so much more
difficult to establish convergence. \\
On the other hand, the {\it stability} property of a finite difference scheme is {\it no}
longer a relational property, but an {\it intrinsic} property which is {\it solely} of the
finite difference scheme itself, therefore, at least in principle, it can be established alone
on the information contained in that finite difference scheme. \\
In this way, in the study of the convergence of finite difference schemes there is clearly a
major interest in establishing a certain connection between the {\it relational} property of
{\it convergence} which is the sought after aim, and on the other hand, the {\it intrinsic}
property of {\it stability}. \\

The above Lax-Richtmyer theorem does establish such a connection, in fact, an equivalence,
between convergence and stability. Unfortunately however, it assumes the {\it completeness} of
the normed space in which all of this happens, in order to be able to prove the implication
"convergent $\Longrightarrow$ stable". \\ \\

{\bf Appendix} \\

We present a simple {\it counterexample} to the celebrated Principle of Uniform Boundedness of
Linear Operators in a Banach Space, based on the fact that the respective normed space fails
to be complete, that is, Banach. This shows that in this principle, the completeness of the
normed space involved is indeed essential. \\

We take the normed space $( X,~||~||~)$ defined as follows \\

(A.1)~~~ $ X ~=~ \{~ x = ( x_0, x_1, x_2, ~.~.~.~ ) \in \mathbb{R}^{\mathbb{N}} ~~\left |~~
                               \begin{array}{l}
                                       \exists~~ m \in \mathbb{N} : \\
                                       \forall~~ n \in \mathbb{N},~ n \geq m : \\
                                       ~~~~ x_n = 0
                                \end{array} ~\right \} $ \\

with the norm given by \\

(A.2)~~~ $ || x || ~=~ \sup~ \{~ | x_n | ~~|~~ n \in \mathbb{N} ~\} $ \\

for $x = ( x_0, x_1, x_2, ~.~.~.~ ) \in X$. \\

Now, for every $k \in \mathbb{N}$, we define the linear operator $T_k : X \longrightarrow X$
by \\

(A.3)~~~ $ T_k ( x_0, x_1, x_2,~.~.~.~ ) ~=~ ( y_0, y_1, y_2,~.~.~.~ ) $ \\

where for $x = ( x_0, x_1, x_2, ~.~.~.~ ) \in X$, we have \\

(A.4)~~~ $ y_n ~=~ \begin{array}{|l} ~~ k x_k ~~~\mbox{if}~~ n = k \\ \\
                                     ~~ 0 ~~~~\mbox{if}~~ n \neq k
                             \end{array} $ \\

It follows easily that \\

(A.5)~~~ $ ||~ T_k ~|| ~=~ k,~~~ k \in \mathbb{N} $ \\

therefore, the family of linear operators $( T_k ~|~ k \in \mathbb{N} )$ is {\it not}
uniformly bounded. \\

On the other hand, given any fixed $x = ( x_0, x_1, x_2, ~.~.~.~ ) \in X$, there exists $m \in
\mathbb{N}$, such that $x_n = 0$, for $n \in \mathbb{N}, n \geq m$. Hence $T_k ( x ) = 0$, for
$k \in \mathbb{N}, k \geq m$. Consequently \\

(A.6)~~~ $ \sup~ \{~ ||~ T_k ( x ) ~|| ~~|~~ k \in \mathbb{N} ~\} ~<~ \infty $ \\

In this way, the family of linear operators $( T_k ~|~ k \in \mathbb{N} )$ is bounded at each
point $x \in X$, and yet it is not uniformly bounded on $X$. \\

The reason for that is obviously in the fact that the normed space $( X,~||~||~)$ is {\it not}
complete. Indeed, $( X,~||~||~)$ is a strict and dense subspace of $l^\infty$. \\

\end{document}